\documentclass{article}
\usepackage{geometry}
\usepackage{amsmath}
\usepackage{amsthm}
\usepackage{amssymb}
\usepackage{mathptmx}
\usepackage{graphicx}
\usepackage{bbding}
\usepackage{tikz}
\theoremstyle{definition}
\newtheorem{defn}{Definition}
\newtheorem*{problem}{Problem}
\newtheorem{example}{Example}

\theoremstyle{remark}
\newtheorem*{rmk}{Remark}

\theoremstyle{plain}

\title{Jenkins-Strebel Differentials on the Riemann Sphere with Four Simple Poles}
\author{Xujia~Chen and Bin~Xu$^\dagger$}
\date{}

\begin{document}
\maketitle
\begin{abstract}
A celebrated and deep theorem in the theory of Riemann surfaces states the existence and uniqueness of the Jenkins-Strebel differentials on a Riemann surface under some conditions, but the proof is non-constructive and examples are difficult to find. This paper deals with an example of a simple case, namely Jenkins-Strebel differentials on the Riemann sphere with four fixed  simple poles. We will give explicit expressions of these Jenkins-Strebel differentials by means of the Weierstrass $\wp$ function and expose a simple algorithm determining the correspondence between these differentials and some classes of simple closed curves on the Riemann sphere with four points removed.
\end{abstract}

\noindent {\bf 2010 Mathematics Subject Classification:}  30F10\\

\noindent {\bf Keywords:} Jenkins-Strebel differential; Weierstrass $\wp$ function; simple closed curve.

\footnote{$^\dagger$
The second author is supported in part by the National Natural Science Foundation of China (Grant Nos. 11571330 and 11271343) and the Fundamental Research Funds for the Central
Universities (Grant No. WK3470000003).}

\section{Introduction}
\label{introduction}
Theorem 25 in an article of Arbarello and Cornalba \cite{arbarello} stated the existence and uniqueness of Jenkins-Strebel differentials with at most simple poles on a Riemann surface under certain conditions. The simplest nontrivial case may be when the surface is the Riemann sphere with four punctures removed. In this case the theorem implies that Jenkins-Strebel differentials on the Riemann sphere with four fixed poles (and holomorphic elsewhere) is in one-to-one correspondence with free homotopy classes of non-oriented closed curves on the surface which contains a simple closed curve separating two of the four punctures apart. However, since proof of the theorem is non-constructive, the explicit correspondence is unknown. In fact, explicit expressions of Jenkins-Strebels arising this way is not easy to find. This paper deals with this simple case and gives explicit expressions of all the Jenkins-Strebel differnentials stated above. Section \ref{JS} gives a brief introduction to Jenkins-Strebel differentials and the result in \cite{arbarello}, section \ref{expression} gives the explicit expressions, and as it will be seen that these Jenkins-Strebel differentials are naturally parametrised by the rational projective line ${\Bbb Q}\cup \{\infty\}$, section \ref{correspondence}  exposes a simple algorithm determining the correspondence between the parameters and the free homotopy classes above.

\section{Jenkins-Strebel differentials}
\label{JS}
We begin by a brief introduction to Jenkins-Strebel differentials. A meromorphic quadratic differential $\omega$ on a Riemann surface $S$ which in local coordinates $(U,\phi,z)$ has form $\omega = f_U (z){dz}^2$ defines a flat metric, the \textit{$\omega$-metric}, which is the metric with local expression $|f|dz\overline{dz}$. Its \textit{horizontal geodesics} are defined to be curves $\gamma=\gamma(t)$ on $S$ satisfying $f_U'(\phi(\gamma(t)))\gamma '(t)^2 > 0$ for every $t$ in $\gamma$'s domain. A \textit{horizontal trajectory} (\textit{trajectory} for short in the following text) is a maximal horizontal geodesic, i.e. is not contained in a longer curve which is also a horizontal geodesic. A trajectory may be a closed curve (\textit{closed}), or connecting two critical points (i.e. zeros and poles) of the quadratic differential (\textit{critical}), or neither. It was shown by Strebel \cite{strebel} that, for the third case (\textit{recurrent}), the closure of the trajectory as a subset of $S$ is of positive measure.

By a \textit{Jenkins-Strebel differential} we refer to a quadratic differential that has no recurrent trajectories. That is, all their trajectories are either closed or critical. The \textit{critical graph} of a Jenkins-Strebel differential is the ribbon graph on the surface consisting of its critical trajectories.

To demonstrate theorem 25 in Arbarello and Cornalba \cite{arbarello}, we need several other concepts first. Both the theorem and the definitions were given in \cite{arbarello}. A \textit{Riemann surface of finite type} is a closed Riemann surface with finitely many punctures. A quadratic meromorphic differential $\omega$ on $S$ is called {\it admissible} if the $\omega$-metric on $S$ has finite area.  An \textit{annular region} is an open region on a Riemann surface isomorphic to $T_R = \{ R\in \mathbb{C} : R<|z|<1 \}$ on the complex plane for some $0 \leq R < 1$. An annular region $\Omega$ is said to be with the same homotopy type as $\gamma$, where $\gamma$ is a closed curve, if $\gamma$ is freely homotopic to a simple closed curve in $\Omega$ which is not homotopically trivial in $\Omega$. It was shown in \cite{arbarello} that if $S$ has negative Euler number, then every closed trajectory of an admissible quadratic differential on a Riemann surface of finite type is contained in a maximal annular region which is swept out by trajectories.

\begin{defn}[Definition 22 in \cite{arbarello}]
Let $S$ be a Riemann surface of finite type. An \textit{admissible system of curves} on $S$ is a collection $(\gamma_1,...,\gamma_k)$ of simple closed curves which are mutually disjoint, neither homotopically trivial in $S$ nor contracting to a puncture of $S$, and such that $\gamma_i$ is not freely homotopic to $\gamma_j$ if $i\neq j$.
\end{defn}

\begin{defn}[Definition 24 in \cite{arbarello}]
Given an admissible system of curves $(\gamma_1,...,\gamma_k)$ on a Riemann surface S of finite type, a collection$(\Omega_1,...,\Omega_k)$ of disjoint subsets of S is said to be a \textit{system of annular regions of type} $(\gamma_1,...,\gamma_k)$ if $\Omega_i$ is either the empty set or an annular region with the same homotopy type as $\gamma_i$, for $i=1,...,k$.
\end{defn}

Now the theorem can be stated:

\newtheorem{thm}{Theorem}
\begin{thm}[Theorem 25 in \cite{arbarello}]
Let $S$ be a Riemann surface of finite type and have negative Euler number. Let $(\gamma_1,...,\gamma_k)$ be an admissible system of curves on $S$ and let $a_1,...a_k$ be positive real numbers. Then there exists on $S$ a unique admissible holomorphic
Jenkins-Strebel differential $\omega$ having the following properties.

i) If $\Gamma$ is the critical graph of $\omega$, then $S\backslash \Gamma = \Omega_1 \cup ... \cup \Omega_k$, where $(\Omega_1,...,\Omega_k)$ is a system of annular regions of type $(\gamma_1,...,\gamma_k)$.

ii) If $\Omega_i$ is not empty, it is swept out by trajectories whose $\omega$-length is $a_i$.
\end{thm}

This manuscript concerns the case when the Riemann surface $S$ is the Riemann sphere with four punctures. More explicitly, consider four distinct points $0,1,\mu,\infty$ on the Riemann sphere, then the surface $S_\mu= \hat{\mathbb{C}} \backslash \{ 0,1,\mu, \infty \}$ satisfies the theorem's condition, and it is obvious that an admissible system of curves on $S_\mu$ can only contain one simple closed curve, which must separates two of $0,1,\mu,\infty$ apart. For later convenience, we make the following definition.

\begin{defn}
By a \textit{pre concerned curve class} we mean a free homotopy class of closed curves on $S_\mu$ that contains a simple closed curve separating two of $0,1,\mu,\infty$ apart.
By a \textit{concerned curve class} we mean an equivalent class of the pre-concerned curve classes module direction. That is, two pre concerned curve classes are  equivalent if and only if one contains $\gamma$ and the other contains $\gamma^{-1}$. When we say a simple closed curve $\gamma$ is \textit{in} a concerned curve class $C$ we just mean that it is in a pre concerned curve class in $C$.
\end{defn}

In addition, since in our case there is only one $a_i$ as in the theorem, and we are not interested in scalar changes of the differential, in the following texts of this section as well as in Section \ref{correspondence}, by Jenkins-Strebel differential we will refer to the equivalence class of Jenkins-Strebel differentials module scalar multiplication (i.e. we may by default assume the closed trajectories have $\omega$-length 1).

For a concerned curve class, choosing a simple closed curve in it as an admissible system of curves on $S_\mu$ we get a Jenkins-Strebel differential, and the four punctures must be simple poles. Different choices of the curve give the same differential, because, for freely homotopic curves (or of opposite direction) $\gamma_1$ and $\gamma_2$ an annular region is of type $\gamma_1$ if and only if of type $\gamma_2$, thus a differential satisfying the theorem's condition for $\gamma_1$ if and only if for $\gamma_2$. Reversely, for a Jenkins-Strebel differential $\omega$ on $\hat{\mathbb{C}}$ with four simple poles at $0,1,\mu,\infty$ and no other poles, a closed trajectory $\gamma$ of $\omega$, as a non-oriented simple closed curve, must separate two of the poles apart. If a simple closed curve on $S_\mu$ is neither homotopically trivial in $S_\mu$ nor contracting to a point in $\{0,1,\mu,\infty\}$, then it intersects $\gamma$ or is freely homotopic to $\gamma$ in $S_\mu$. On the other hand, since any two closed trajectories do not intersect, every closed trajectory of $\omega$ must be freely homotopic to $\gamma$.
Thus $\omega$ gives a unique concerned curve class, which obviously will also give $\omega$ in the process above. Therefore, Theorem 1 tells us that {\it Jenkins-Strebel differentials on $\hat{\mathbb{C}}$ with four simple poles at $0,1,\mu,\infty$ are in one-to-one correspondence with the concerned curve classes}. Here comes our problem to be dealt with in this paper:

\begin{problem}
Find explicit expressions for the Jenkins-Strebel differentials on the Riemann sphere with four simple poles at $0,1,\mu,\infty$ and no other poles.
\end{problem}

\paragraph{}
This problem is solved in Section \ref{expression}, by means of pulling back to complex tori.

\section{Expression}
\label{expression}
A quadratic differential $\omega$ on the Riemann sphere with four simple poles can be written in the form $\omega=\frac{a}{z(z-1)(z-\mu)} dz^2$. Our only problem is to determine that for which $a$ it is Jenkins-Strebel, that is, has no recurrent horizontal trajectory. To do this, we pull $\omega$ back through a branched double cover from the torus to $\hat{\mathbb{C}}$ that has four branch values $0,1,\mu,\infty$. The simple poles being branch values of degree two, their pre-images are regular points for the pulled-back differential. Pre-images of regular points remain regular for the pulled-back differential. Therefore the pulled-back differential has no critical points.

The double cover we need is the Weierstrass $\wp$ function. A \textit{lattice} on the complex plane is a subset \{$m\omega_1+n\omega_2$ : $m,n$ are integers\} for some non-zero complex numbers $\omega_1,\omega_2$ with $\omega_1/\omega_2 \notin \mathbb{R}$. (Differently chosen basis $\omega_1,\omega_2$ may differ by a unimodular transformation.) For a fixed lattice the related Weierstrass $\wp$ function is a meromorphic function on $\mathbb{C}$ with its period being the lattice points, poles being lattice points, and branching points being half lattice points (with degree 2). (Thus it is actually defined on the torus.) The branch values are $e_1:=\wp(\omega_1/2),e_2:=\wp(\omega_2/2),e_3:=\wp((\omega_1+\omega_2)/2)$ and $\infty$. The modular lambda function is defined to be $\lambda(\omega_1/\omega_2)=(e_3-e_2)/(e_1-e_2)$.
Take basis change into consideration, $\lambda=\lambda(\tau)$ is actually defined on $(\mathbb{C} \backslash \mathbb{R}) / G$, where $G$ is the congruence subgroup mod 2 of the modular group (i.e. $\lambda'=\frac{a\lambda+b}{c\lambda+d}$ with $a,b$ odd integers, $b,c$ even integers, $ad-bc=1$). We also know that $\lambda=\lambda(\tau)$ is holomorphic on ${\Bbb C}\backslash {\Bbb R}$ and could achieve any complex numbers other than $0$ and $1$. Reference: \cite{ahlfors}.

For the problem, consider a lattice with $\lambda(\tau)=\mu$. For convenience we normalise it to assume that it is spanned by $1$ and $\tau$, where $\tau$ is in $\lambda^{-1}(\mu)$. Since the $\wp$ function has branch values $e_1,e_2,e_3,\infty$, let us first consider the quadratic differentials with simple poles there. Such a quadratic differential can be expressed as $\omega=\frac{c}{(z-e_1)(z-e_2)(z-e_3)} dz^2$. Pull it back through $\wp$ and compute using the following property of the $\wp$ function
\begin{equation}
\wp(z)'^2=4(\wp(z)-e_1)(\wp(z)-e_2)(\wp(z)-e_3)
\end{equation}
to get
\begin{align}
\wp^\ast \omega&=\frac{c}{(\wp(w)-e_1)(\wp(w)-e_2)(\wp(w)-e_3)}{\wp(w)'}^2 dw^2 \\
&=4c \ dw^2
\end{align}
where $w$ is the natural coordinate (up to transformation) on the torus as quotient of the complex plane by the lattice. Therefore the horizontal trajectories $\gamma$ are determined by
\begin{equation}
4c \ \gamma'(w)^2 > 0.
\end{equation}
$\omega$ has no recurrent horizontal trajectory if and only if $\wp^\ast\omega$ has none, and the trajectory $\gamma$ is not recurrent if and only if $arg(\gamma')=arg(1+q\tau)$ for some rational number $q$ or $arg(\gamma')=arg(\tau)$, so $\omega$ is Jenkins-Strebel if and only if either there is a rational $q$ and a positive real $k$ such that $c=k/(1+q\tau)^2$ or there is a positive real $k$ such that $c=k/\tau^2$, as illustrated in Figure 1. 

~\\
\begin{figure}[h]
\begin{centering}
\begin{tikzpicture}
\pgftransformcm{1}{0}{0.4}{1}{\pgfpoint{0cm}{0cm}}
\draw[style=help lines,dashed] (-1,-1) grid[step=1cm] (2,2);
\foreach \x in {-1,0,1,2}{
	\foreach \y in {-1,0,1,2}{
		\node [draw,circle,inner sep=1pt,fill] at (\x,\y) {};
	}
}
\node [below] at (0,0) {0};
\node [below] at (1,0) {1};
\node [left] at (0,1) {$\tau$};
\node [below right] at (1,0.7) {$1+q\tau$};
\node [draw,circle,inner sep=1pt] at (1,0.7) {};
\draw (-1.25,-0.875) -- (2.5,1.75);
\node [below] at (2.5,1.75) {$\gamma$};
\end{tikzpicture}
\caption{Trajectory on torus}
\end{centering}
\label{fig1}
\end{figure}
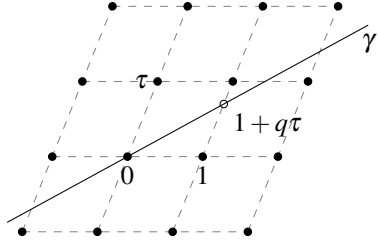

Now for the case presented in the problem. The Mobius transformation $T(z)=(z-e_2)/(e_1-e_2)$ sends $e_1, e_2, e_3, \infty$ to 1, 0, $\mu$, $\infty$ respectively, so the quadratic differential $\omega=\frac{a}{z(z-1)(z-\mu)} dz^2$ which has four simple poles $0,1,\mu,\infty$ is pulled back through $T$ to be
\begin{align}
T^\ast \omega &= \frac{a}{T(z)(T(z)-1)(T(z)-\mu)} T'(z)^2 dz^2 \\
&= \frac{a}{\frac{z-e_2}{e_1-e_2} (\frac{z-e_2}{e_1-e_2}-1) (\frac{z-e_2}{e_1-e_2}-\mu)} \frac{1}{(e_1-e_2)^2} dz^2 \\
&= \frac{a(e_1-e_2)}{(z-e_1)(z-e_2)(z-e_3)}dz^2,
\end{align}
which has four simple poles $e_1, e_2, e_3, \infty$.
Therefore $\omega$ is Jenkins-Strebel if and only if either $a=\frac{k}{(1+q\tau)^2 (e_1-e_2)}$ for some rational $q$ and real positive $k$ or $a=\frac{k}{\tau^2(e_1-e_2)}$ for some real positive $k$, where $\tau$ is an element in $\lambda^{-1}(\mu)$, and $e_1, e_2$ are determined by $\tau$ as above. For the latter case we just say this is when $q$ is $\infty$ for convenience.

\section{Correspondence}
\label{correspondence}
Recall that in this section by Jenkins-Strebel differential we will mean equivalent class of Jenkins-Strebel differentials module scalar change. We have already seen at the end of Section \ref{JS} that on $S_\mu$ the Jenkins-Strebel differentials we study are in one-to-one correspondence with the concerned curve classes. As seen in Section \ref{expression}, fix a $\tau$ in $\lambda^{-1}(\mu)$ and the differentials can be parametrised by $\mathbb{Q} \cup \{ \infty \}$. Thus so do the concerned curve classes. This parameterisation is faithful, because as seen in the formula in Chaper \ref{expression}, different $q$ give different differentials. This can also be seen by pulling back to the torus: different $q$ correspond to closed curves not freely homotopic on the torus. What's more, free homotopy classes of curves are in one-to-one correspondence with conjugacy classes of the fundamental group. We conclude that for a fixed $\tau$ the following four are in one-to-one correspondence: $\mathbb{Q} \cup \{ \infty \}$, Jenkins-Strebel differentials, the concerned curve classes, conjugacy classes of $\pi_1(S_\mu)$ module inverse which contains a simple closed curve separating two of $0,1,\mu,\infty$ apart. In what follows in this section we will give a simple algorithm determining the conjugacy class of $\pi_1(S_\mu)$ for each $q \in \mathbb{Q} \cup \{ \infty \}$.

\begin{figure}
\centering
\begin{tikzpicture}
\pgftransformcm{1}{0}{0.4}{1}{\pgfpoint{0cm}{0cm}}
\draw[style=help lines,dashed] (-1,-1) grid[step=0.5cm] (2,2);

\foreach \x in {-1,0,1}{
	\foreach \y in {-0.5,0.5,1.5}{
		\fill[gray] (\x,\y) rectangle (\x+0.5,\y+0.5);
		\draw [red,thick] (\x,\y) -- (\x+0.5,\y);
		\draw [blue,thick] (\x+0.5,\y) -- (\x+0.5,\y+0.5);
		\draw [green,thick] (\x+0.5,\y+0.5) -- (\x,\y+0.5);
		\draw [yellow,thick] (\x,\y+0.5) -- (\x,\y);
	}
}
\foreach \x in {-0.5,0.5,1.5}{
	\foreach \y in {-1,0,1}{w
		\fill[gray] (\x,\y) rectangle (\x+0.5,\y+0.5);
		\draw [green,thick] (\x,\y) -- (\x+0.5,\y);
		\draw [yellow,thick] (\x+0.5,\y) -- (\x+0.5,\y+0.5);
		\draw [red,thick] (\x+0.5,\y+0.5) -- (\x,\y+0.5);
		\draw [blue,thick] (\x,\y+0.5) -- (\x,\y);
	}
}

\foreach \x in {-1,0,1,2}{
	\foreach \y in {-1,0,1,2}{
		\node [draw,circle,inner sep=1pt,fill] at (\x,\y) {};
	}
}
\node [below] at (0,0) {0};
\node [below] at (1,0) {1};
\node [left] at (0,1) {$\tau$};
\draw (-1.1,-0.8) -- (2.35,1.9);
\node [below] at (2.35,1.9) {$\gamma$};
\end{tikzpicture}
\caption{Chessboardwise-colored lattice}
\label{lattice}
\end{figure}
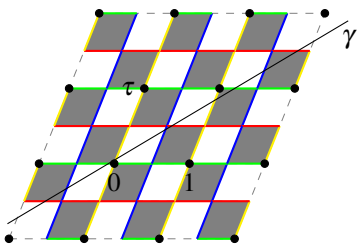

\begin{figure}
\centering
\begin{tikzpicture}
\draw [green,thick] (-1,-1) to [out=15,in=165] (1,-1);
\draw [blue,thick] (1,-1) to [out=105,in=255] (1,1);
\draw (1,1) to [out=195,in=-15] (-1,1);
\draw (-1,1) to [out=-75,in=75] (-1,-1);
\draw [red,thick] (1,1) to [out=205,in=-25] (-1,1);
\draw [yellow,thick] (-1,1) to [out=-65,in=65] (-1,-1);
\fill [gray] (1,1) to [out=195,in=-15] (-1,1) to [out=-25,in=205] (1,1);
\fill [gray] (-1,1) to [out=-75,in=75] (-1,-1) to [out=65,in=-65] (-1,1);
\node [left] at (-1,1) {0};
\node [left] at (-1,-1) {$\infty$};
\node [right] at (1,-1) {1};
\node [right] at (1,1) {$\mu$};
\foreach \x in {-1,1}{
	\foreach \y in {-1,1}{
		\node [draw,circle,inner sep=1.5pt,fill] at (\x,\y) {};
	}
}

\node [draw,circle,inner sep=0.5pt,fill] at (0,0) {};
\node [below right] at (0,0) {$x_0$};

\draw [->] (0,0) to [out=25,in=-95] (0.885,0.5);
\draw [dashed] (0.885,0.5) to [out=95,in=-5] (0.5,0.885);
\draw (0.5,0.885) to [out=190,in=65] (0,0);
\node [right] at (0.885,0.5) {b};

\draw [->] (0,0) to [out=115,in=-5] (-0.5,0.885);
\draw [dashed] (-0.5,0.885) to [out=185,in=85] (-0.885,0.5);
\draw (-0.885,0.5) to [out=280,in=155] (0,0);
\node [above] at (-0.5,0.885) {c};

\draw [->] (0,0) to [out=205,in=85] (-0.885,-0.5);
\draw [dashed] (-0.885,-0.5) to [out=275,in=175] (-0.5,-0.885);w
\draw (-0.5,-0.885) to [out=10,in=245] (0,0);
\node [left] at (-0.885,-0.5) {$d=c^{-1}b^{-1}a^{-1}$};

\draw [->] (0,0) to [out=295,in=175] (0.5,-0.885);
\draw [dashed] (0.5,-0.885) to [out=5,in=265] (0.885,-0.5);
\draw (0.885,-0.5) to [out=100,in=-25] (0,0);
\node [below] at (0.5,-0.885) {a};
\end{tikzpicture}
\caption{$\hat{\mathbb{C}} \backslash \{0,1,\mu,\infty\}$ as a pillowcase}
\label{pillowcase}
\end{figure}

Color the complex plane chessboardwise as shown in Figure \ref{lattice}. For a white or gray tile, $\wp$ is a bijection when restricted to it, thus a homeomorphism when restricted to the interior of it. $\wp$ is also continuous on the closed tile, so the image of the tile under $\wp$ is a topological square with four vertices branch values. What's more, due to the periodicity and the even property of $\wp$, tiles of the same color have the same image rigeon. Therefore, $\hat{\mathbb{C}}$ is topologically a pillowcase with the two faces being images of tiles of the two colors respectively and the four corners being branch values. Denote $T=(z-e_2)/(e_1-e_2)$ the Mobius transformation the same as in Section \ref{expression}. Composing $\wp$ with $T$ does not change anything stated above, and the branching values become $0,1,\mu,\infty$, as shown in Figure \ref{pillowcase}. Coloring of the edges shows the image-preimage relation. From then on we will just call $\hat{\mathbb{C}} \backslash \{0,1,\mu,\infty\}$ \textit{pillowcase} and call the image of white tiles the \textit{front face}, image of gray tiles the \textit{back face} of the pillowcase. The loops $a,b,c$ represent a basis of $\pi_1(\hat{\mathbb{C}} \backslash \{0,1,\mu,\infty\}, x_0)$ for $x_0$ a point on the front face.

Now we can give the algorithm. For a given $q$ in $\mathbb{Q} \cup \{ \infty \}$, arbitrarily choose a closed horizontal trajectory of the corresponding Jenkins-Strebel differential and denote the pre-image of it under $T \circ \wp$ as $\gamma$. (Such curves are straight lines of slope $q$ on the latticed plane, and we will call them \textit{$q$-lines}.) Choose an arbitrary point $y_0$ on $\gamma$ which is also in a white tile, then go along $\gamma$ (in arbitrary direction) from $y_0$ until the image of $\gamma$ under $T\circ \wp$ goes back to the image of $y_0$ (assume it is $x_0$, without loss of generality). $\gamma$ goes through white and gray tiles alternatively as its image goes through the front and back face of the pillowcase. Since we are focusing on the homotopy property of the image of $\gamma$, it does not matter if we assume that it comes back to $x_0$ every time passing through the front face. Therefore the image of $\gamma$ can be divided into several loops starting from $x_0$, passing through the back face once and going back to $x_0$. The pre-image of such a loop (up to homotopy) is a segment on $\gamma$ starting from a white tile, passing through a gray tile, ending in the next white tile. It is easy to see that the homotopy class of such a loop is determined by the two edges it passes, as shown in Table \ref{edgecorrs}. Therefore, we just need to write down the colors of the edges $\gamma$ passed across in order, pair them sequently, find the elements in the fundamental group corresponding to the pairs, and multiply them in order to get a fundamental group element whose conjugate class (module inverse) is what we want.

\begin{table}
\centering
\caption{Passing edges -- fundamental group elements correspondence}
\begin{tabular}{cccccc}
\hline
red$\to$yellow & yellow$\to$green & green$\to$blue & blue$\to$red & red$\to$green & yellow$\to$blue \\
$c$ & $d$ & $a$ & $b$ & $cd$ & $da$ \\
\hline
yellow$\to$red & green$\to$yellow & blue$\to$green & red$\to$blue & green$\to$red & blue$\to$yellow \\
$c^{-1}$ & $d^{-1}$ & $a^{-1}$ & $b^{-1}$ & $d^{-1}c^{-1}$ & $a^{-1}d^{-1}$ \\
\hline
\end{tabular}
\label{edgecorrs}
\end{table}

\begin{example}[$q=0$]
When $q=0$, $\gamma$ is a horizontal line on the complex plane, and $T \circ \wp (\gamma)$ goes to the back face of the pillowcase only once. The edges passed across are blue to yellow, corresponding to $a^{-1}d^{-1}=bc$.
\end{example}

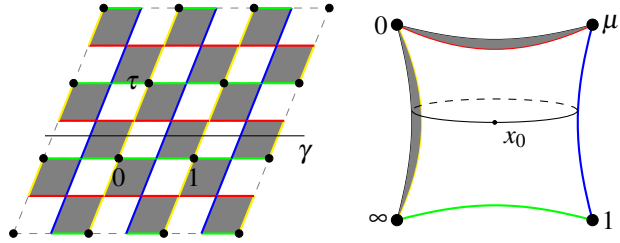
\begin{figure}[h]
\centering

\begin{tikzpicture}
\pgftransformcm{1}{0}{0.4}{1}{\pgfpoint{0cm}{0cm}}
\draw[style=help lines,dashed] (-1,-1) grid[step=0.5cm] (2,2);

\foreach \x in {-1,0,1}{
	\foreach \y in {-0.5,0.5,1.5}{
		\fill[gray] (\x,\y) rectangle (\x+0.5,\y+0.5);
		\draw [red,thick] (\x,\y) -- (\x+0.5,\y);
		\draw [blue,thick] (\x+0.5,\y) -- (\x+0.5,\y+0.5);
		\draw [green,thick] (\x+0.5,\y+0.5) -- (\x,\y+0.5);
		\draw [yellow,thick] (\x,\y+0.5) -- (\x,\y);
	}
}
\foreach \x in {-0.5,0.5,1.5}{
	\foreach \y in {-1,0,1}{w
		\fill[gray] (\x,\y) rectangle (\x+0.5,\y+0.5);
		\draw [green,thick] (\x,\y) -- (\x+0.5,\y);
		\draw [yellow,thick] (\x+0.5,\y) -- (\x+0.5,\y+0.5);
		\draw [red,thick] (\x+0.5,\y+0.5) -- (\x,\y+0.5);
		\draw [blue,thick] (\x,\y+0.5) -- (\x,\y);
	}
}

\foreach \x in {-1,0,1,2}{
	\foreach \y in {-1,0,1,2}{
		\node [draw,circle,inner sep=1pt,fill] at (\x,\y) {};
	}
}
\node [below] at (0,0) {0};
\node [below] at (1,0) {1};
\node [left] at (0,1) {$\tau$};
\draw (-1.1,0.3) -- (2.35,0.3);
\node [below] at (2.35,0.3) {$\gamma$};
\end{tikzpicture}
\ \ \
\begin{tikzpicture}[scale=1.3]
\draw [green,thick] (-1,-1) to [out=15,in=165] (1,-1);
\draw [blue,thick] (1,-1) to [out=105,in=255] (1,1);
\draw (1,1) to [out=195,in=-15] (-1,1);
\draw (-1,1) to [out=-75,in=75] (-1,-1);
\draw [red,thick] (1,1) to [out=205,in=-25] (-1,1);
\draw [yellow,thick] (-1,1) to [out=-65,in=65] (-1,-1);
\fill [gray] (1,1) to [out=195,in=-15] (-1,1) to [out=-25,in=205] (1,1);
\fill [gray] (-1,1) to [out=-75,in=75] (-1,-1) to [out=65,in=-65] (-1,1);
\node [left] at (-1,1) {0};
\node [left] at (-1,-1) {$\infty$};
\node [right] at (1,-1) {1};
\node [right] at (1,1) {$\mu$};
\foreach \x in {-1,1}{
	\foreach \y in {-1,1}{
		\node [draw,circle,inner sep=1.5pt,fill] at (\x,\y) {};
	}
}

\node [draw,circle,inner sep=0.5pt,fill] at (0,0) {};
\node [below right] at (0,0) {$x_0$};

\draw (-0.845,0.12) arc (180:360:0.845 and 0.12);
\draw [dashed] (0.845,0.12) arc (0:180:0.845 and 0.12);


\end{tikzpicture}
\caption{The $q=0$ example}
\end{figure}

\begin{figure}[h]
\centering

\begin{tikzpicture}
\pgftransformcm{1}{0}{0.4}{1}{\pgfpoint{0cm}{0cm}}
\draw[style=help lines,dashed] (-1,-1) grid[step=0.5cm] (2,2);

\foreach \x in {-1,0,1}{
	\foreach \y in {-0.5,0.5,1.5}{
		\fill[gray] (\x,\y) rectangle (\x+0.5,\y+0.5);
		\draw [red,thick] (\x,\y) -- (\x+0.5,\y);
		\draw [blue,thick] (\x+0.5,\y) -- (\x+0.5,\y+0.5);
		\draw [green,thick] (\x+0.5,\y+0.5) -- (\x,\y+0.5);
		\draw [yellow,thick] (\x,\y+0.5) -- (\x,\y);
	}
}
\foreach \x in {-0.5,0.5,1.5}{
	\foreach \y in {-1,0,1}{w
		\fill[gray] (\x,\y) rectangle (\x+0.5,\y+0.5);
		\draw [green,thick] (\x,\y) -- (\x+0.5,\y);
		\draw [yellow,thick] (\x+0.5,\y) -- (\x+0.5,\y+0.5);
		\draw [red,thick] (\x+0.5,\y+0.5) -- (\x,\y+0.5);
		\draw [blue,thick] (\x,\y+0.5) -- (\x,\y);
	}
}

\foreach \x in {-1,0,1,2}{
	\foreach \y in {-1,0,1,2}{
		\node [draw,circle,inner sep=1pt,fill] at (\x,\y) {};
	}
}
\node [below] at (0,0) {0};
\node [below] at (1,0) {1};
\node [left] at (0,1) {$\tau$};
\draw [dashed] (-0.6,-1.2) -- (1.5,3);
\draw (-0.5,-1.3) -- (1.5,2.7);
\node [below] at (-0.5,-1.3) {$\gamma$};
\node [draw,circle,inner sep=0.8pt,fill] at (0.25,0.2) {};
\node [draw,circle,inner sep=0.8pt,fill] at (1.25,2.2) {};
\node [below] at (0.25,0.2) {$y_0$};
\node [below] at (1.25,2.2) {$y_0'$};
\end{tikzpicture}
\ \ \
\begin{tikzpicture}[scale=1.3]
\draw [green,thick] (-1,-1) to [out=15,in=165] (1,-1);
\draw [blue,thick] (1,-1) to [out=105,in=255] (1,1);
\draw (1,1) to [out=195,in=-15] (-1,1);
\draw (-1,1) to [out=-75,in=75] (-1,-1);
\draw [red,thick] (1,1) to [out=205,in=-25] (-1,1);
\draw [yellow,thick] (-1,1) to [out=-65,in=65] (-1,-1);
\fill [gray] (1,1) to [out=195,in=-15] (-1,1) to [out=-25,in=205] (1,1);
\fill [gray] (-1,1) to [out=-75,in=75] (-1,-1) to [out=65,in=-65] (-1,1);
\node [left] at (-1,1) {0};
\node [left] at (-1,-1) {$\infty$};
\node [right] at (1,-1) {1};
\node [right] at (1,1) {$\mu$};
\foreach \x in {-1,1}{
	\foreach \y in {-1,1}{
		\node [draw,circle,inner sep=1.5pt,fill] at (\x,\y) {};
	}
}

\node [draw,circle,inner sep=0.5pt,fill] at (0,0) {};
\node [below right] at (0,0) {$x_0$};
\draw (-0.6,-0.9) to [out=15,in=195] (0.6,0.9);
\draw [dashed] (0.6,0.9) to [out=-15,in=90] (0.885,0.5);
\draw (0.885,0.5) to [out=-105,in=15] (0.2,-0.85);
\draw [dashed] (0.2,-0.85) to [out=165,in=-15] (-0.6,0.9);
\draw (-0.6,0.9) to [out=-170,in=90] (-0.843,0);
\draw [dashed] (-0.843,0) to [out=-90,in=-180] (-0.6,-0.9);

\end{tikzpicture}
\caption{The $q=2$ example}
\label{q=2}
\end{figure}
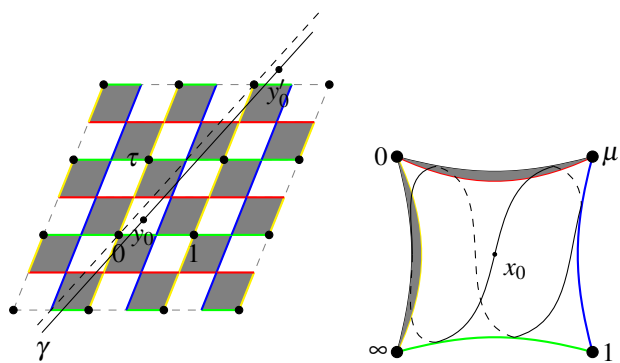

\begin{example}[$q=2$]
When $q=2$, as shown in Figure \ref{q=2}, the image of $\gamma$ passes through the edges in $red \to blue \to green \to red \to yellow \to green$, which corresponds to $b^{-1} (d^{-1}c^{-1}) d =babc^{-1}b^{-1}a^{-1}$.
\end{example}

\begin{rmk} Since a closed horizontal trajectory (when given a direction) as a simple closed curve separates two of $0,1,\mu,\infty$ apart,
three cases can occur concerning which two points are in the same region. Therefore $\mathbb{Q} \cup \{\infty\}$ may also be divided into three parts correspondingly. For a $q$ in $\mathbb{Q} \cup \{\infty\}$, write it as $m/n$ where $m,n$ are coprime integers. (For $q=\infty$ make it $1/0$, for $q=0$ make it $0/1$.) The $q-line$ passing through the origin on the latticed plane is mapped by $T \circ \wp$ to a critical horizontal trajectory of the Jenkins-Strebel differential corresponding to $q$, which connects $\infty$, the image of the origin, to another critical point of the differential. Since this two critical points are connected by a critical horizontal trajectory which cannot intersect a closed horizontal trajectory, they must lie in the same region enclosed by a closed horizontal trajectory. Therefore $\mathbb{Q} \cup \{\infty\}$ is divided into three parts determined by that, of $0,1,\mu$, whose preimage under $T \circ \wp$ the $q-line$ through origin also passes through, but this is easy to be seen determined by the parity of $m,n$. They cannot be simultaneously even, and the other three cases correspond to the three parts of  $\mathbb{Q} \cup \{\infty\}$.
\end{rmk}

\bibliographystyle{plain}

~\\

\noindent Xujia Chen\\
School for the Gifted Young,
University of Science and Technology of China\\
No. 96 Jinzhai Road, Hefei,
Anhui Province\  230026\  P. R. China\\
chenxj94@mail.ustc.edu.cn\\
\noindent {\it Current address}: Department of Mathematics,
Stony Brook University,
Stony Brook, NY 11974\\
xujia.chen@stonybrook.edu\\

\noindent Bin Xu\\
Wu Wen-Tsun Key Laboratory of Mathematics, Chinese Academy of Sciences\\
School of Mathematical Sciences,
University of Science and Technology of China\\
No. 96 Jinzhai Road, Hefei,
Anhui Province\  230026\  P. R. China\\
\Envelope \,bxu@ustc.edu.cn

\end{document}